\documentclass[12pt]{amsart} 
\usepackage[table]{xcolor}
\definecolor{light-gray}{gray}{0.85}

\usepackage{amsfonts,graphics,amsmath,amsthm,amsfonts,amscd, amssymb,amsmath,latexsym,bm}
\usepackage{epsfig}
\usepackage{flafter}
\usepackage{float}

\usepackage{mathtools}

\usepackage{lipsum}

\usepackage[english,italian]{babel}
\usepackage{graphicx,subfigure}
\usepackage[a4paper,top=30mm,bottom=30mm,left=30mm,right=30mm]{geometry}
\usepackage{tikz}
\usetikzlibrary{matrix,arrows,positioning,decorations.pathmorphing}
\tikzset{
    >=stealth',
    punkt/.style={
           rectangle,
           rounded corners,
           draw=black, thick,
           text width=3em,
           minimum height=2em,
           text centered},
    pil/.style={
           ->,
           thick,
           shorten <=2pt,
           shorten >=2pt,}
}

\usepackage[pdftex,
            pdfauthor={Stefano Filipazzi Franco Rota},
            pdftitle={An example of Berglund-H\"ubsch mirror symmetry for a Calabi-Yau complete intersection},
            pdfsubject={Algebraic geometry, mirror symmetry, Landau-Ginzburg mirror symmetry},
            pdfkeywords={filipazzi rota  mirror symmetry berglund hubsch krawitz calabi yau}]{hyperref}

\usepackage[babel]{csquotes}
\usepackage[maxnames=99,style=alphabetic,backend=bibtex]{biblatex}
\usepackage{colortbl}

\usepackage[T1]{fontenc}
\usepackage{longtable}
\usepackage{cases}
\usepackage{booktabs}
\usepackage{caption}
\usepackage{vmargin}
\usepackage{verbatim}

\theoremstyle{plain}
\newtheorem{theorem}{Theorem}[section]

\newtheorem{definition-lemma}[theorem]{Definition-Lemma}
\newtheorem{question}[theorem]{Question}
\newtheorem{defn}[theorem]{Definition}

\theoremstyle{definition}

\newtheorem*{ack*}{Acknowledgements}
\newtheorem{remark}[theorem]{Remark}
\newtheorem{example}[theorem]{Example}

\newcommand{\C}{\mathbb{C}}			
\newcommand{\Z}{\mathbb{Z}}			
\newcommand{\rar}{\rightarrow}		

\DeclareMathOperator{\Jac}{Jac}		
\DeclareMathOperator{\spn}{span}		

\def\O#1.{\mathcal {O}_{#1}}			
\def\af #1.{\mathbb A^{#1}}				
\def\ses#1.#2.#3.{0\to #1\to #2\to #3 \to 0}		
\def\xrar#1.{\xrightarrow{#1}}			
\def\K#1.{K_{#1}}						
\def\bM#1.#2.{\mathbf{M}_{#1,#2}}				
\def\bL#1.#2.{\mathbf{L}_{#1,#2}}				
\def\bB#1.#2.{\mathbf{B}_{#1,#2}}				
\def\subs#1.{_{#1}}						
\def\sups#1.{^{#1}}						

\newcommand{\Addresses}{{
  \bigskip
  \footnotesize

  S.~Filipazzi, \textsc{Department of Mathematics, University of Utah,
    Salt Lake City,\\ UT 84112, USA}\par\nopagebreak
  \textit{E-mail address}: \texttt{filipazz@math.utah.edu}
  
  F.~Rota, \textsc{Department of Mathematics, University of Utah,
    Salt Lake City,\\ UT 84112, USA}\par\nopagebreak
  \textit{E-mail address}: \texttt{rota@math.utah.edu}
}}

\newcommand{\pr}[1]{\mathbb{P}^{#1}}
\newcommand{\SL}{\text{SL}}

\title{An example of Berglund-H\"ubsch mirror symmetry for a Calabi-Yau complete intersection}
\author{Stefano Filipazzi, Franco Rota}
\date{\today}
\thanks{The two authors were partially supported by NSF FRG Grant DMS-1265285.}

\begin{document}

\selectlanguage{english}

\begin{abstract}
We study an example of complete intersection Calabi-Yau threefold due to Libgober and Teitelbaum \cite{LT93}, and verify mirror symmetry at a cohomological level. Direct computations allow us to propose an analogue to the Berglund-H\"ubsch mirror symmetry setup for this example \cite{BH93}. We then follow the approach of Krawitz to propose an explicit mirror map \cite{Kra10}.
\end{abstract}

\maketitle
\tableofcontents

\section{Introduction}

This note is an account of a calculation which we carried out during the Summer School ``Pragmatic 2015'', and represents a first step in our project of generalizing the Krawitz-Chiodo-Ruan cohomological isomorphisms to Berglund-H\"ubsch mirror pairs of Calabi-Yau complete intersections.

In 1993, Berglund and H\"ubsch describe a procedure to construct the mirror manifold of a large class of Calabi-Yau hypersurfaces in weighted projective space, using the formalism of Landau-Ginzburg models \cite{BH93}. Later, Krawitz gives an explicit description of the mirror map between the state spaces of two Berglund-H\"ubsch mirror LG models \cite{Kra10}, and Borisov re-proves and generalizes his statements \cite{Bor13}.

Since the work by Berglund and H\"ubsch, the role of Landau-Ginzburg models in mirror symmetry has increased in importance, particularly after Fan, Jarvis and Ruan defined their quantum invariants in \cite{FJR13} (FJRW theory), and extended the definition to an even more general setting in \cite{FJR15}. In a subsequent work, Chiodo and Nagel recognize the state space of Landau-Ginzburg isolated singularities defined in \cite{FJR13} as an instance of relative orbifold Chen-Ruan cohomology \cite{CN15}. This allows to define the cohomology of hybrid Landau-Ginzburg models, and to identify it with the cohomology of Calabi-Yau complete intersections in weighted projective space. In this sense, \cite{CN15} extends the LG/CY correspondence, which was proven to hold for hypersurfaces in \cite{CR11}, to the case of complete intersection Calabi-Yau. The situation in the hypersurface case can be summarized by the following diagram,
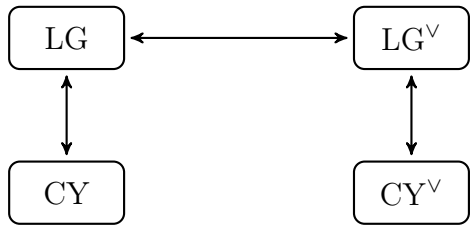
\begin{center}
\begin{tikzpicture}[node distance=1cm, auto,]
 \node[punkt] (CY) {CY};
 \node[punkt, right=3cm of CY] (CYv) {$\mathrm{CY}^\vee$};
 \node[punkt, above=1.2cm of CY] (LG) {LG}
 edge[pil, <->] (CY);
 \node[punkt, above=1.2cm of CYv] (LGv) {$\mathrm{LG}^\vee$}
 edge[pil, <->] (LG);
 \path (LGv) edge[pil, <->] (CYv);
\end{tikzpicture}
\end{center}
where an arrow indicates an isomorphism of the corresponding state spaces. Existence of the vertical arrows has been established in \cite{CR11}, while the horizontal one is the object of \cite{Kra10} and \cite{Bor13}.

In the case of a complete intersection, \cite{CN15} provides the vertical arrows of the diagram
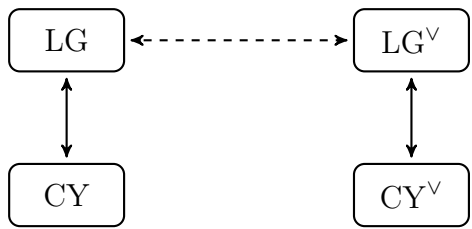
\begin{center}
\begin{tikzpicture}[node distance=1cm, auto,]
 \node[punkt] (CY) {CY};
 \node[punkt, right=3cm of CY] (CYv) {$\mathrm{CY}^\vee$};
 \node[punkt, above=1.2cm of CY] (LG) {LG}
 edge[pil, <->] (CY);
 \node[punkt, above=1.2cm of CYv] (LGv) {$\mathrm{LG}^\vee$}
 edge[pil, <->, dashed] (LG);
 \path (LGv) edge[pil, <->] (CYv);
\end{tikzpicture}
\end{center}
by proving an orbifold version of the Thom isomorphism. A related result by Libgober \cite{Lib15} shows the invariance of the elliptic genus for complete intersections and other examples of GIT quotients. 

The object of this work is to exhibit a horizontal arrow in a particular example by Libgober and Teitelbaum \cite{LT93}. In this way, we aim at making progress towards a definition of a mirror symmetry construction for Calabi-Yau complete intersections of the same type as the Berglund-H\"ubsch construction. We describe a complete intersection Calabi-Yau threefold and its mirror orbifold using an analogue of the Berglund-H\"ubsch formalism, and propose a generalization of Krawitz's mirror map. 

The direction of this paper is suggested by the classical result of Batyrev and Borisov on cohomological mirror symmetry for Calabi-Yau complete intersections \cite{BB94}, although at the moment it is unclear how to explicitly relate to their work.

Starting from a paper by Libgober and Teitelbaum \cite{LT93}, we consider the complete intersection of two cubics in $\pr 5$ together with its mirror. We recall the setup for hybrid Landau-Ginzburg models in Section \ref{section LG models}. In Section \ref{section description example}, using the Chiodo-Nagel CY/LG-correspondence \cite{CN15}, we rephrase the example of \cite{LT93} in terms of hybrid Landau-Ginzburg models. Then, in Section \ref{section Krawitz revisited}, we generalize the construction of Krawitz for the quintic Fermat hypersurface in $\pr 4$ with the language proposed for complete intersections. Finally, in Section \ref{section computations}, we verify that the state spaces of the complete intersection and of its mirror are isomorphic, and present an explicit mirror map generalizing Krawitz's formula.

\begin{ack*}
We started working on this project during the PRAGMATIC 2015 Research School in Algebraic Geometry and Commutative Algebra, ``Moduli of Curves and Line Bundles'', held in Catania, Italy, in July 2015. We are very grateful to Alfio Ragusa, Francesco Russo, and Giuseppe Zappal\`a, the organizers of the PRAGMATIC school, for the wonderful environment they created at the school. We thank the lecturers of the school, Alessandro Chiodo and Filippo Viviani. In particular, we are indebted with Alessandro Chiodo for first suggesting this problem, and then advising us through it. We are also grateful to J\'{e}r\'{e}my Gu\'{e}r\'{e} for helpful insight and comments.
\end{ack*}

\section{Hybrid Landau-Ginzburg models} \label{section LG models}

In this section, we recall the notion of hybrid Landau-Ginzburg models. Our setup follows closely the one in \cite{CN15}, whose notation we adopt for the most part. The main difference is that we consider a broader class of groups of symmetries. This is in line with the groups considered in \cite{LT93}.

Let $\bm W = \left\lbrace W_1,\ldots ,W_r \right\rbrace$ be a set of quasi-homogeneous polynomials of degrees $d_1,\ldots ,d_r$ defining a complete intersection $$X_{\bm W}\coloneqq \left\lbrace W_1=\ldots =W_r=0 \right\rbrace \subset \mathbb P(w_1,\ldots ,w_n)$$ in weighted projective space. Assume furthermore that $X_{\bm W}$ is {\it non-degenerate}, i.e. 
\begin{itemize}
\item the choice of weights $w_i$ is unique;
\item $X_{\bm W}$ is smooth outside the origin.
\end{itemize} The \emph{maximal group of diagonal symmetries} of $X_{\bm W}$ is the maximal subgroup $\Gamma_{max}\subset (\C^*)^{n}$ of diagonal matrices $\gamma$ preserving $X_{\bm W}$. More explicitly, $\gamma$ has diagonal entries of the form $(\gamma_1,\ldots ,\gamma_n)=\alpha\bar{\lambda}\coloneqq (\alpha_1\lambda^{w^1},\ldots , \alpha_n\lambda^{w_n})$, where $\lambda\in\C^*$, and $\alpha=\alpha\bar{1}$ is a vector whose entries are all non-zero. By definition, for any $1 \leq i \leq r$ and admissible choice of $\alpha$, there exists $\beta_{i,\alpha}\in \C^*$ such that $W_i(\alpha_1x_1,\ldots ,\alpha_nx_n)= \beta_{i,\alpha} W_i(x_1,\ldots ,x_n)$. By a \emph{group of diagonal symmetries} we mean a subgroup of $\Gamma_{max}$.

\begin{example} \label{group J}
The group $\Gamma_{max}$ contains the following special element. Let $d$ denote the greatest common divisor of the $d_i$, and define
\begin{equation*}
J \coloneqq  \begin{pmatrix}
    e^{2\pi i \frac{w_1}{d}}  &  \\
   & \ddots & \\
   && e^{2\pi i \frac{w_n}{d}}
\end{pmatrix}.
\end{equation*}
The element $J$ is analogous to the generator of the group $\bm \mu_d$ in the hypersurface case. The group $\left\langle J \right\rangle$ will play an important role later.  
\end{example}

\begin{example}
Among the subgroups of $\Gamma_{max}$, one can consider the special linear group of diagonal symmetries of a complete intersection, denoted $\SL (\bm W)$. It contains elements of $\Gamma_{max}$ of determinant 1. 
\end{example}

\begin{defn}\label{notation}
Let $\gamma$ be a diagonal symmetry acting on $\C^n$ as a diagonal matrix with entries $\exp(2\pi i \gamma_i)$, with $\gamma_i\in \left[0,1\right)$. Then, we denote $\gamma$ by $(\gamma_1,\ldots ,\gamma_n)$. If we have $\gamma_i=\frac{c_i}{d}$ for all $i$, we write $(\gamma_1,\ldots ,\gamma_n)=\frac{1}{d}(c_1,\ldots,c_n)$. The \emph{age} of the element $\gamma$, denoted $a_\gamma$, is the sum $\sum_i \gamma_i$. 
\end{defn}

Unless otherwise stated, we will follow the classic literature (see \cite[Ch.~3]{CR14}), and assume that any group of diagonal symmetries $\Gamma$ contains the torus 
$$\Gamma_0 \coloneqq \left\lbrace \bar\lambda=(\lambda^{w_1},\ldots ,\lambda^{w_n}) | \lambda\in\C^* \right\rbrace .$$
Then, $\Gamma_0$ is the connected component of the identity in $\Gamma$, and we denote by $G$ the quotient $\Gamma/\Gamma_0$.

Now, consider $\C \sups n+r.$ with coordinates $(x_1, \ldots, x_n,p_1, \ldots, p_r)$. We can extend the action of $\Gamma$ to $\C^{n+r}$ by setting
\begin{equation}\label{extendaction}
 \alpha\bar{\lambda}(\bm x, \bm p) \coloneqq (\alpha_1\lambda^{w^1}x_1,\ldots , \alpha_n\lambda^{w_n}x_n, \beta_{1,\alpha}^{-1}\lambda^{-d_1}p_1,\ldots ,\beta_{r,\alpha}^{-1}\lambda^{-d_r}p_r).
\end{equation}
It is useful to introduce a notation for the fixed points of an element $\gamma\in\Gamma$. We set
\begin{equation*}
\C^n_\gamma\coloneqq \left\lbrace \bm x\in\C^n| \gamma \cdot \bm x = \bm x \right\rbrace, \qquad n_\gamma\coloneqq \text{dim}\C^n_\gamma,
\end{equation*}
where the action of $\gamma$ is restricted to the first set of coordinates. Similarly, we define
\begin{equation*}
\C^r_\gamma\coloneqq \left\lbrace \bm p\in\C^r| \gamma \cdot \bm p = \bm p \right\rbrace, \qquad r_\gamma\coloneqq \text{dim}\C^r_\gamma.
\end{equation*}
The integer $n_\gamma$ (respectively $r_\gamma$) counts the dimension of the space spanned by the $\bm x$ (resp. $\bm p$) variables that are fixed by $\gamma$. For a polynomial $V\in \C[\bm x, \bm p]$, we let $V_\gamma$ denote $V|_{\C^n_\gamma\times\C^r_\gamma}$.

For an element $V\in \C[\bm x, \bm p]$, its {\it chiral algebra} $\mathcal{Q}_V$ is defined as
\begin{equation*}
\mathcal Q_{V} \coloneqq d\bm x\wedge d \bm p \otimes \Jac(V),
\end{equation*}
where we are formally tensoring the {\it Jacobi ring} 
$$
\Jac(V) \coloneqq \dfrac{\C[x_1,\ldots,x_n,p_1,\ldots,p_r]}{(\frac{\partial V}{\partial x_1},\ldots,\frac{\partial V}{\partial x_n},\frac{\partial V}{\partial p_1},\ldots,\frac{\partial V}{\partial p_r})}
$$
by the top form $d\bm x \wedge d \bm p \coloneqq dx_1\wedge \ldots \wedge dx_n \wedge dp_1\wedge \ldots \wedge dp_r$. For brevity, in the following we will often omit the wedge symbol and write $d\bm x d \bm p \coloneqq d\bm x \wedge d \bm p $. We assign bidegree $(D-k,k)$ to the elements of $\mathcal Q_V$ which have degree $k$ in the $\bm p$ variables, where $D=n-r-1$, .

By construction, $\gamma^* (p_iW_i)=p_iW_i$ for all $\gamma\in\Gamma$ and all $i$.
Then we can define the $\Gamma$-invariant function
\begin{equation*}
\begin{split}
\overline{W} : \C^{n+r}  &\rar \C\\
  (\bm{x},\bm{p})  &\mapsto p_1W_1(\bm x)+\ldots +p_rW_r(\bm x).
\end{split}
\end{equation*}
Let $M \coloneqq \left\lbrace \overline W = t_0 \right\rbrace$ denote the fiber over any point $t_0\neq 0$.
Consider the open set in $\C^{n+r}$ defined as $U_{LG} \coloneqq \C^n\times (\C^r\smallsetminus \left\lbrace\bm 0\right\rbrace) $. The quotient of $U_{LG}$ by the action of $\C^*$ is the total space of the vector bundle
\[ \left[ U_{LG}/\C^* \right] = \bigoplus\limits_{j=1}^n \mathcal O_{\pr{}(\bm{d})}(-w_i)\eqqcolon \mathcal O_{\bm d}(-\bm w). \]

\begin{defn}The datum $(\mathcal O_{\bm d}(-\bm w), \overline{W}, \Gamma)$ is the hybrid Landau-Ginzburg model with superpotential
\begin{equation*}
\overline{W}:\left[ \mathcal O_{\bm d}(-\bm w)/G \right] \to \C.
\end{equation*}
The generalized state space of the hybrid Landau-Ginzburg model $(\mathcal O_{\bm d}(-\bm w), \overline{W})$ is 
\[ \mathcal H_\Gamma^{p,q}(W_1,\ldots ,W_r)\coloneqq  H_{\mathrm{CR}}^{p+r,q+r}\left( \left[ \mathcal O_{\bm d}(-\bm w)/G \right], \left[ F/G \right] \right), \]
where $F$ is the quotient stack $\left[ M/\C^* \right]$ inside $\mathcal O_{\bm d}(-\bm w)$.

\end{defn}

For any two non-negative integers $r\geq n$ we need to consider the bigraded ring $dt(-n)\otimes \C[t]/(t^{r-n})$, where the element $dt(-n)\otimes t^k$ has bidegree $(n+k,n+k)$. We will make use of the following theorem due to Chiodo and Nagel \cite{CN15}. 

\begin{theorem}[{\cite[Theorem 4.3]{CN15}}]\label{thm_CN}
Let $W_i$, $r$ and $n$ be as introduced above. Then, we have 
\[ \mathcal H_\Gamma^{*}(W_1,\ldots ,W_r)=\bigoplus_{\gamma\in \Gamma}H_\gamma(-a_\gamma+r), \]
where $H_\gamma$ with its double grading is given by 
   \begin{subnumcases}{H_\gamma=}
   \left(\mathcal Q_{\widetilde W_\gamma}\right)^\Gamma & $\mbox{ if } r_\gamma<n_\gamma$ \label{Jacobi sector}
   \\
  dt(-n_\gamma)\otimes \C[t]/(t^{r_\gamma-n_\gamma}) & $\mbox{ if } r_\gamma \geq n_\gamma$. \label{projective sector}
\end{subnumcases}
and carries the Tate twist $(-a_\gamma + r)$. 
\end{theorem}

\section{Description of the example} \label{section description example}

In this section, we introduce the complete intersection Calabi-Yau object of this work. This example was first described by Libgober and Teitelbaum \cite{LT93}.

Consider $\pr5$ with coordinates $[x_1:x_2:x_3:X_1:X_2:X_3]$. Let $X\subset \pr5$ be the Calabi-Yau threefold defined by the vanishing of the polynomials
\begin{equation*}
\begin{split}
W_{1} &= x_1^3 + x_2^3 + x_3^3 - 3X_1X_2X_3\\
W_{2} &= X_1^3 + X_2^3 + X_3^3 - 3x_1x_2x_3.
\end{split}
\end{equation*}
Since we are in the standard projective space $\pr5$, and we are considering two polynomials of degree 3, the group $\langle J \rangle$ introduced in Example \ref{group J} is isomorphic to $\bm{\mu}_3$.

Introduce two new variables, $p_1$ and $p_2$, and form the polynomial $W \coloneqq p_1W_1 + p_2W_2$. Extend the action of the group $\left\langle J \right\rangle$ to $\C^6\times\C^2$ by letting $J$ act with weights $(1,\ldots,1,-3,-3)$. Similarly, following equation (\ref{extendaction}), we can extend the action of any group $\Gamma$ of diagonal symmetries to the variables $p_1$ and $p_2$.

According to the hypersurface case of Berglund-H\"ubsch mirror symmetry (see \cite[Ch. 3]{CR14}), to study the generalized state spaces of the Hybrid Landau-Ginzburg model induced by $W_1$ and $W_2$, we have to make choices for a mirror set of polynomials and for a pair of groups. As suggested by Libgober and Teitelbaum, we choose the same set of polynomials for the mirror. Indeed, this is analogous to the Fermat hypersurface case. Then, let $\Gamma \coloneqq \Gamma_0\cdot \left\langle J \right\rangle$, and $\Gamma^T \coloneqq \Gamma_{max}\cdot \SL(\bm W)$. This choice is in line with the Berglund-H\"ubsch prescription (see Rmk. \ref{G81}).

By Theorem \ref{thm_CN}, the only elements of a group of diagonal symmetries $\Sigma$ contributing to the generalized state space of $(\mathcal O_{\bm d}(-\bm w), \overline{W}, \Sigma)$ are the elements $\sigma$ fixing at least one of the coordinates $(x_1,x_2,x_3,X_1,X_2,X_3,p_1,p_2)$.

By direct inspection, one checks that the only elements of $\Gamma$ that fix at least one of the variables are the elements of $\left\langle J \right\rangle$. Analogously, when considering $\Gamma^T$, one is left with those elements $\gamma \in \SL(\bm W)$ that fix at least one of the variables. These $\gamma$ generate a subgroup $\Gamma'\subset \SL(\bm W)$ whose elements $g_{\lambda,\mu,\alpha,\beta,\delta,\epsilon}$ act on $(x_1,x_2,x_3,X_1,X_2,X_3,p_1,p_2)$ via multiplication by
\[(\lambda \xi_3^\alpha\xi_9^\mu, \lambda \xi_3^\beta\xi_9^\mu, \lambda \xi_9^\mu, \lambda \xi_3^{-\delta}\xi_9^{-\mu},\lambda \xi_3^{-\epsilon}\xi_9^{-\mu},\lambda \xi_9^{-\mu},\lambda^{-3}\xi_3^{-\mu},\lambda^{-3}\xi_3^{\mu}).\] 
Here $\xi_3$ is a primitive third root of unity, $\xi_9$ is a primitive ninth root of unity, and $\lambda\in\bm\mu_9$. Furthermore, the condition 
$$\alpha + \beta = 3\mu = \delta+ \epsilon \mod 3$$
is satisfied.

\begin{remark}\label{G81}
The group $\Gamma^T$ also contains $\Gamma_0$, and the quotient $\Gamma^T/\Gamma_0$ coincides with the group $G_{81}$ described in \cite{LT93}. We have $\Gamma/\Gamma_0=\left\langle J \right\rangle$. 
\end{remark}

As a consequence of Libgober and Teitelbaum's work  \cite{LT93} and an application of Chiodo and Nagel's CY/LG-correspondence \cite[Thm. 5.1]{CN15}, one knows that the generalized state spaces $\mathcal H_{\Gamma}^{p,q}(W_1,W_2)$ and $ \mathcal H_{\Gamma^T}^{3-p,q}(W_1,W_2)$ have the same dimension. The object of this work is to further investigate this correspondence, and propose a possible mirror map. The main idea is summarized in the following.

\begin{defn} \label{main result}
We define an explicit mirror map $\mathcal H_{\Gamma^T}^{1,1}(W_1,W_2) \to \mathcal H_{\Gamma}^{2,1}(W_1,W_2)$, generalizing Krawitz's map for hypersurfaces.
\end{defn}

\begin{remark}
The mirror map in Definition \ref{main result} is expected to extend to a map $\mathcal H_{\Gamma^T}^{*}(W_1,W_2) \to \mathcal H_{\Gamma}^{*}(W_1,W_2)$. At this point, the main difficulty is to choose appropriate representatives among the 73 generators of $\mathcal H_{\Gamma}^{1,2}(W_1,W_2)$ in order to recognize a pattern similar to the one occurring among the generators of $\mathcal H_{\Gamma}^{2,1}(W_1,W_2)$, which are listed in Table \ref{tab_2,1}.
\end{remark}

\section{Krawitz's mirror map revisited} \label{section Krawitz revisited}

In this section, we discuss how the map mentioned in Definition \ref{main result} is related to the work of Krawitz. More precisely, we show how the approach suggested for Calabi-Yau complete intersections is a generalization of the one already known for hypersurfaces.

We will focus on a concrete example, the Fermat quintic $X_W \subset \pr4$.  This is given as the vanishing locus of the polynomial
$$
W=x_1^5+x_2^5+x_3^5+x_4^5+x_5^5.
$$
First, we consider the group $\langle J \rangle$ as group of diagonal symmetries. Thus, we first consider the pair $(W,\langle J \rangle)$. Berglund-H\"ubsch mirror symmetry guarantees the existence of a mirror model, denoted by $(W^T,\langle J \rangle^T)$ \cite{Kra10, CR14}. In our case, $W=W^T$, and $\langle J \rangle^T= \SL(W^T)= \SL (W)$ is a well known group, denoted by $G_{125}$.

Now, we would like to regard the hypersurface $X_W$ as a ``complete intersection of one hypersurface''. Therefore, we introduce one auxiliary variable $p$, and, as explained in equation (\ref{extendaction}), extend the action of the diagonal symmetries acting on $X_W$ to $p$. We can then regard $pW$ as a function on $\C^5 \times \C^*$, and consider the induced Landau-Ginzburg model.

Among the groups fixing $X_W$, there is $\C^*\cdot \left\langle J \right\rangle$. In particular, $\C^*$ acts on $x_i$ with weight 1, and on $p$ with weight $-5$. We can compute the corresponding state spaces. By Theorem \ref{thm_CN}, only elements $\gamma$ fixing some variable contribute to the computations. Furthermore, it is easy to check that, in this case, $\gamma \in \C^*$ is forced to satisfy $|\gamma|=1$. More precisely, $\gamma$ is either 1, or a primitive fifth-root of unity. In the first case, all the variables are fixed, while, in the latter one, just $p$ remains unchanged.

The element $\gamma=1$ contributes to the state space with a summand of the form
$$
\left(d\bm p d \bm x \otimes \frac{\C[p,x_1,\ldots,x_5]}{(px_1^4,\ldots,px_5^4,\sum \subs i=1.^5 x_i^5)} \right)  \sups \C^*..
$$
Its elements carry a bigrading $(3-k,k)$, where $k$ is the degree of the variable $p$ in the element. As $d\bm p d \bm x$ is $\C^*$ invariant, and we are looking for $\C^*$-invariant elements, we need polynomial coefficients containing five among $x_1, \ldots , x_5$ for every appearance of $p$. There is just a one-dimensional family in degree $(3,0)$, there are spaces with dimension 101 in degree $(2,1)$ and $(1,2)$, and there is one more one-dimensional family in degree $(0,3)$.

If $\gamma$ is a primitive fifth-root of unity, we have $r_\gamma=1$ and $n_\gamma =0$. Therefore, each one of these elements contributes with a one-dimensional vector space, generated by one element, denoted by $\left. 1|g \right\rangle$.

Taking into account the grading of the above pieces, we recover the well known Hodge diamond

\begin{center}
\begin{tikzpicture}[yscale = .7, xscale = .7]
\node at (0,2) {1} ;
\node at (0,1) {1} ;
\node at (0,-1) {1} ;
\node at (0,-2) {1} ;
\node at (1,0) {101} ;
\node at (2,0) {1} ;
\node at (-1,0) {101} ;
\node at (-2,0) {1} ;
\end{tikzpicture}
\end{center}

Now, we want to consider the mirror of $(pW,\C^*\cdot \left\langle J \right\rangle)$. As $W$ is a polynomial of Fermat type, we will let the mirror polynomial be $pW$ as well. As $\C^*$ is the minimial group of diagonal symmetries containing $\C^*$, we consider the maximal group of symmetries fixing the polynomial $pW$. Let $\Gamma$ be such group. It can be checked by direct computation that $\Gamma= \C^* \cdot G_{125}$. It acts as
$$
(x_1,x_2,x_3,x_4,x_5,p) \mapsto (\lambda x_1, \lambda \xi^a x_2, \lambda \xi^b x_3, \lambda \xi^c x_4, \lambda \xi \sups -a-b-c. x_5, \lambda \sups -5. p),
$$
where $\lambda \in \C^*$, $\xi$ is a primitive fifth-root of unity, and $a,b,c \in \Z$.

\begin{remark}
In \cite{CN15}, Chiodo and Nagel consider the maximal group fixing $W$, and then they extend its action to $pW$. On the other hand, for our purposes, such a group turns out to be too rigid, and does not recover the right cohomologies.
\end{remark}

As above, the only elements $\gamma \in \Gamma$ contributing to the state spaces are the ones where $\lambda$ is a fifth-root of unity. In particular, $p$ is fixed by any relevant element $\gamma \in \Gamma$. This highlights how the variable $p$ is not necessary for the computations in the hypersurface case.

Since $\lambda$ needs to be a fifth-root of unity, we can absorb it in the contribution of $\xi$. Therefore, any element $\gamma$ contributing to the computations can be encoded by a 6-tuple
$$
(a_1,a_2,a_3,a_4,a_5;0),
$$
where the action on $x_i$ is given by multiplication by $\xi \sups a_i.$. We will use this formalism as a slight modification of the one in Definition \ref{notation}, in order to distinguish the $p$ and the $x$ variables.

Now, we are left with direct computations. We will proceed by cases.
\begin{itemize}
\item Assume that $p$ is the only variable fixed. In such case, $r_\gamma = 1$, and $n_\gamma = 0$. Therefore, each element contributes with a one-dimensional vector space. There are 204 of these elements, and their ages groups them in four families of dimension 1, 101, 101 and 1 respectively. These compute the twisted sector of the cohomology ring.

\item If $p$ and some but not all of the $x_i$ are fixed, there is no contribution.

\item The identity element, represented by the string $(0,0,0,0,0;0)$, contributes to the untwisted sector. We have to consider the ring
$$
\left(d\bm p d \bm x \otimes \frac{\C[p,x_1,\ldots,x_5]}{(px_1^4,\ldots,px_5^4,\sum \subs i=1.^5 x_i^5)} \right)  \sups \C^*\cdot G_{125}..
$$
This ring has dimension four as vector space, with basis
$$
(d\bm p d \bm x\otimes 1, d\bm p d \bm x\otimes p x_1 \ldots x_5, d\bm p d \bm x\otimes (p x_1 \ldots x_5)^2, d\bm p d \bm x\otimes (p x_1 \ldots x_5)^3).
$$
Each one of these elements has a different bidegree.
\end{itemize}

If we put the above facts together, we recover the Hodge diamond

\begin{center}
\begin{tikzpicture}[yscale = .7, xscale = .7]
\node at (0,2) {1} ;
\node at (0,1) {101} ;
\node at (0,-1) {101} ;
\node at (0,-2) {1} ;
\node at (1,0) {1} ;
\node at (2,0) {1} ;
\node at (-1,0) {1} ;
\node at (-2,0) {1} ;
\end{tikzpicture}
\end{center}

With the elements explicitly listed, we can proceed and rewrite Krawitz's mirror map with this new notation. First, we consider the generators of the untwisted sector of the mirror. These are mapped according to
$$
(px_1 \ldots x_5)^{a-1} \otimes d \bm p d \bm x \left| id \right\rangle \mapsto 1 \left| (a,a,a,a,a;0) \right\rangle,
$$
where $a \in \lbrace 1,2,3,4 \rbrace$.

Symmetrically, the twisted elements of the mirror are mapped to the untwisted sector of the quintic. The map is given by 
\begin{equation} \label{Krawitz eqtn}
1 \left| (a,b,c,d,e;0) \right\rangle \mapsto d \bm p d \bm x \otimes p \sups \frac{a+b+c+d+e}{5}-1. x_1 \sups a-1. x_2 \sups b-1. x_3 \sups c-1. x_4 \sups d-1. x_5 \sups e-1.\left| id \right\rangle.
\end{equation}

\begin{remark}
The quantity $\frac{a+b+c+d+e}{5}$ is the age $a_\gamma$ of the group element $\gamma$ considered. Therefore, under the mirror map introduced by Krawitz \cite{Kra10}, the auxiliary variable $p$ appears with exponent $a_\gamma -1$.
\end{remark}

\section{The mirror map} \label{section computations}

In this section, we construct the mirror map in Definition \ref{main result}. As we noticed in the example of the Fermat quintic, there is no harm in dropping the assumption that the diagonal symmetries have determinant 1. Indeed, the new group elements coming into the picture do not contribute to the computation of the state spaces. Recall that we consider the groups $\Gamma=\Gamma_0\cdot \left\langle J \right\rangle$ and $\Gamma^T$ acting on the same set of polynomials (see Section \ref{section description example}).

In the following, we first compute the state space associated to $\Gamma$. Then, we consider the one associated to $\Gamma^T$, and we conclude describing the mirror map.

\subsection{The state spaces $\mathcal H_{\Gamma}^{p,q}(W_1,W_2)$}\label{sect_Gamma}

We follow the notation of Theorem \ref{thm_CN}. Notice that for $\gamma \in \Gamma$, if $r_\gamma=n_\gamma=0$, then $H_\gamma=0$ gives no contribution to the state space. Then, only the three elements of $\left\langle J \right\rangle \subset \Gamma$ give non-trivial contributions. We will proceed by cases.
\begin{itemize}
\item The non-trivial elements have ages 2 and 4 respectively. They  both contribute with the cohomology of a projective line $\pr{}(3,3) \cong \pr{1}$, and form the twisted sectors of $\mathcal H_{\Gamma}^{*}(W_1,W_2)$.

\item The contribution of the identity element is encoded in $\mathcal Q_{W}^{\Gamma}$. This ring carries a bigrading $(D-k,k)$, where $D=n-r-1=3$, and $k \in \lbrace 0,1,2,3 \rbrace$ is the degree in the variables $p_1$ and $p_2$. The relations in $\text{Jac}(V)$ are:
\begin{align*}
x_1^3+x_2^3+x_3^3- 3 X_1X_2X_3,\\
3 x_1x_2x_3 - X_1^3+X_2^3+X_3^3,\\
p_1x_i^2 - p_2x_jx_k,\\
p_2X_i^2 - p_1X_jX_k,
\end{align*}
where all the $i,j,k$ are assumed to be distinct. The form 
$$
d\bm p d\bm x d \bm X \coloneqq dp_1 \wedge dp_2 \wedge dx_1 \wedge dx_2 \wedge dx_3 \wedge dX_1 \wedge dX_2 \wedge dX_3$$
is invariant under the action of $\Gamma$. Hence, the degree $(3,0)$ component is generated by $d\bm p d\bm x d \bm X \otimes 1$. The 73 independent generators of $\mathcal H_{\Gamma}^{2,1}(W_1,W_2)$ are listed in Table \ref{tab_2,1}. They determine the whole untwisted sector, which is completed by $\mathcal H_{\Gamma}^{1,2}(W_1,W_2)$ and $\mathcal H_{\Gamma}^{0,3}(W_1,W_2)$, whose dimensions are 73 and 1 respectively.

\end{itemize}

We can summarize the above computations with the following Hodge diamond.

\begin{center}
\begin{tikzpicture}[yscale = .7, xscale = .7]
\node at (0,2) {1} ;
\node at (0,1) {1} ;
\node at (0,-1) {1} ;
\node at (0,-2) {1} ;
\node at (1,0) {73} ;
\node at (2,0) {1} ;
\node at (-1,0) {73} ;
\node at (-2,0) {1} ;
\end{tikzpicture}
\end{center}

\subsection{The state spaces $\mathcal H_{\Gamma^T}^{p,q}(W_1,W_2)$}\label{sect_Gamma^vee}

Now, we consider the mirror side. We will list the elements of $\Gamma^T$, and study their contributions to the state spaces. As explained in Definition \ref{notation}, we will write $(\gamma_1,\ldots ,\gamma_8)$ if $\gamma$ acts on $\C^6 \times \C^2$ as a diagonal matrix with entries $\exp(2\pi i \gamma_i)$.

As showed in Theorem \ref{thm_CN}, the contribution of an element $\gamma$ has two different behaviors, depending on how $r_\gamma$ and $n_\gamma$ compare to each other. Therefore, we will proceed by cases.

\begin{itemize}
\item If $r_\gamma < n_\gamma$, we are in the case given by equation \eqref{Jacobi sector}. These elements are listed in Table $\ref{tab_JacobiGamma^vee}$. One of these elements is the identity of $\Gamma^T$. It corresponds to the untwisted sector of $\mathcal H_{\Gamma^T}^{*}(W_1,W_2)$, which is 4-dimensional and generated by powers of $d\bm p d\bm x d \bm X \otimes p_1x_1^3$ in the Jacobi ring.

The other four elements all give a two dimensional Jacobi ring, and behave in a similar way. For instance, $\gamma=\frac19(0,0,0,3,3,3;0,0)$ induces a direct summand $H_\gamma=\spn (d\bm p d\bm x d \bm X \otimes x_1^3,d\bm p d\bm x d \bm X \otimes x_2^3)$, which contributes as two-dimensional subspace of $\mathcal H_{\Gamma^T}^{1,1}(W_1,W_2)$. A second one of these elements contributes with another two-dimensional subspace of $\mathcal H_{\Gamma^T}^{1,1}(W_1,W_2)$. Similarly, the two remaining elements correspond to two two-dimensional subspaces of $\mathcal H_{\Gamma^T}^{2,2}(W_1,W_2)$.

\item If $r_\gamma \geq n_\gamma$, we are in the case given by equation \eqref{projective sector}. Each element in Table $\ref{tab_projGamma^vee}$ contributes to the state space with $r_\gamma - n_\gamma$ elements. If $r_\gamma - n_\gamma=2$, the two generators belong one to 
$$\mathcal H_{\Gamma^T}^{a_\gamma-r+n_\gamma,a_\gamma-r+n_\gamma}(W_1,W_2),
$$
and the other to
$$
\mathcal H_{\Gamma^T}^{a_\gamma+1-r+n_\gamma,a_\gamma+1-r+n_\gamma}(W_1,W_2).
$$
There are four such elements. The first one of them contributes with one generator to $\mathcal H_{\Gamma^T}^{0,0}(W_1,W_2)$, and with one to $\mathcal H_{\Gamma^T}^{1,1}(W_1,W_2)$. Similarly, a second element generates the top cohomology, contributing with one generator to $\mathcal H_{\Gamma^T}^{3,3}(W_1,W_2)$, and with one to $\mathcal H_{\Gamma^T}^{2,2}(W_1,W_2)$. Then, the two remaining elements behave in the same way: Each one of them has one generator lying in $\mathcal H_{\Gamma^T}^{1,1}(W_1,W_2)$, and one in $\mathcal H_{\Gamma^T}^{2,2}(W_1,W_2)$.

On the other hand, if $r_\gamma - n_\gamma=1$, there is only one generator, which belongs to 
$$\mathcal H_{\Gamma^T}^{a_\gamma-r+n_\gamma,a_\gamma-r+n_\gamma}(W_1,W_2).
$$
The last column in Table \ref{tab_projGamma^vee} counts the elements of the {\it same type}. For example, from the element $\frac19(3,3,3,3,0,6;0,0)$ we can obtain 5 more by permuting $X_1,X_2$ and $X_3$ via the action of the symmetric group $S_3$. These new elements are considered of the same type as the starting one, and are not explicitly listed in Table \ref{tab_projGamma^vee}. They are taken into account by counting 6 elements, including their representative, of the type of $\frac19(3,3,3,3,0,6;0,0)$.
\end{itemize}


We can summarize the above computations with the following Hodge diamond.

\begin{center}
\begin{tikzpicture}[yscale = .7, xscale = .7]
\node at (0,2) {1} ;
\node at (0,1) {73} ;
\node at (0,-1) {73} ;
\node at (0,-2) {1} ;
\node at (1,0) {1} ;
\node at (2,0) {1} ;
\node at (-1,0) {1} ;
\node at (-2,0) {1} ;
\end{tikzpicture}
\end{center}

\subsection{The mirror map}

We present the explicit isomorphism between the vector spaces $\mathcal H_{\Gamma^T}^{1,1}(W_1,W_2)$ and $\mathcal H_{\Gamma}^{2,1}(W_1,W_2)$ mentioned in Definition \ref{main result}. Its main feature is that it respects the structure emerging from the lists in Tables \ref{tab_2,1}, \ref{tab_JacobiGamma^vee} and \ref{tab_projGamma^vee}, where it becomes evident that these two vector spaces admit a further interesting subdivisions into direct summands. 

We follow the notation in \cite{Kra10} and Theorem \ref{thm_CN}: Given an element $\sigma$ of a group of diagonal symmetries $\Sigma$, we denote by $\omega|\sigma\rangle$ the element of the state space $\mathcal H^* _\Sigma$ determined by $\omega \in H_\sigma$. 

Now, we explain the mirror map via its action on some elements listed in Table \ref{tab_projGamma^vee}, which represent the different possible behaviors. For the complete action of the mirror map, we refer to Table \ref{tab_complete_mirrormap} in the Appendix.

First, we notice that there is a correspondence between direct summands in $\mathcal H_{\Gamma^T}^{1,1}(W_1,W_2)$ and $\mathcal H_{\Gamma}^{2,1}(W_1,W_2)$. For the reader's convenience, we will list them relying on Tables \ref{tab_2,1}, \ref{tab_JacobiGamma^vee} and \ref{tab_projGamma^vee}.
\begin{itemize}
\item In $\mathcal H_{\Gamma}^{2,1}(W_1,W_2)$, we isolated six groups of generators, each one consisting of nine elements. These correspond to the six groups of nine elements $\gamma \in \Gamma^T$ having $r_\gamma =1$, $n_\gamma = 0$, and $a_\gamma = 3$.

\item The two groups of six generators of $\mathcal H_{\Gamma}^{2,1}(W_1,W_2)$ correspond to the two groups of six elements $\gamma \in \Gamma^T$ having $r_\gamma =2$, $n_\gamma = 1$, and $a_\gamma = 2$.

\item In $\Gamma^T$ there are two elements with $r_\gamma=2$, $n_\gamma = 3$ and $a_\gamma = 1$. These contribute with two two-dimensional vector spaces. On the side of $\mathcal H_{\Gamma}^{2,1}(W_1,W_2)$, these correspond to the vector spaces $\langle p_1 X_1^3 , p_1 X_2^3 \rangle$ and $\langle p_2 x_1^3 , p_2 x_2^3 \rangle$. Notice that, by the relations in the Jacobi ring, these spaces contain $p_1 X_3^3$ and $p_2 x_3^3$ respectively.

\item To conclude, among the elements $\gamma \in \Gamma^T$ with $r_\gamma =2$ and $n_\gamma =0$, three of these contribute with a one-dimensional subspace of $\mathcal H_{\Gamma^T}^{1,1}(W_1,W_2)$ each. These will correspond to the three sporadic polynomials listed at the top of Table \ref{tab_2,1}.
\end{itemize}

Before getting into the details of the mirror map, it is worth asking the following interesting question.

\begin{question}
Do the aforementioned direct summands of the spaces $\mathcal H_{\Gamma^T}^{1,1}(W_1,W_2)$ and $\mathcal H_{\Gamma}^{2,1}(W_1,W_2)$ have an interpretation from the point of view of geometry or physics?
\end{question}

Now, we will proceed to describing the mirror map, analyzing each one of the cases listed above.

The two first groups of elements present a fundamental common feature, namely $r_\gamma - n_\gamma = 1$. This is essentially the only case occurring in the twisted sector for the hypersurface case \cite{Kra10}. As showed in equation \eqref{Krawitz eqtn}, once we introduce a $p$ in the hypersurface case, the exponents of the $x_i$ contribute to the exponent of this new variable. We will follow this idea, bearing in mind that $p_1$ is related to the set of variables $x_i$, while $p_2$ to the set $X_i$.

For the elements with $r_\gamma - n_\gamma = 1$, we then propose the assignment
\begin{equation} \label{assignment rule}
\begin{split}
dt\left|\frac{1}{9}(b_1,b_2,b_3,c_1,c_2, c_3;a_1,a_2)\right\rangle &\mapsto \\
d \bm p d \bm x d \bm X \otimes p_1 \sups \min \lbrace \lfloor \frac{b_j}{3} \rfloor \rbrace. p_2 \sups \min \lbrace \lfloor \frac{c_j}{3} \rfloor \rbrace.  \prod_{i=1}^3 x_i \sups \lfloor \frac{b_i}{3} \rfloor -\min \lbrace \lfloor \frac{b_j}{3} \rfloor \rbrace. & X_i \sups \lfloor \frac{c_i}{3} \rfloor -\min \lbrace \lfloor \frac{c_j}{3} \rfloor \rbrace.  \left|id \right\rangle .
\end{split}
\end{equation}
Indeed, from Table \ref{tab_projGamma^vee}, it is clear that exactly one set of variables has an action which is a rotation of at least $e \sups \frac{2\pi i}{3}.$ on each single coordinate. We record this information choosing the corresponding variable $p_i$ to appear in the assigned polynomial. Then, in order to choose the exponent appearing for each variable in the polynomial in output, we compare the action withing each group $x_i$ and $X_i$. This is where the normalizing exponents $-\min \lbrace \lfloor \frac{b_j}{3} \rfloor \rbrace$ and $-\min \lbrace \lfloor \frac{c_j}{3} \rfloor \rbrace$ come from. Finally, we point out that looking at quantities of the form $\lfloor \frac{n}{3}\rfloor$ is consistent with distinguishing the action of third-roots of unity from the one of primitive ninth-roots.

As example of the assignments described in equation \eqref{assignment rule}, we have
$$dt\left|\frac{1}{9}(3,3,3,3,0,6;0,0)\right\rangle \mapsto d \bm p d \bm x d \bm X \otimes p_1X_1X_3^2 |id\rangle ,$$
and
$$dt\left|\frac{1}{9}(2,2,5,6,3,6;3,0)\right\rangle \mapsto d \bm p d \bm x d \bm X \otimes p_2x_3X_1X_3 |id\rangle.$$

Now, we can study how the mirror map acts on the subspaces corresponding to the elements $\gamma \in \Gamma^T$ with $r_\gamma =2$, $n_\gamma =3$, and $a_\gamma =1$. The first element is
$$
\frac{1}{9}(0,0,0,3,3,3;0,0),
$$
which contributes with the subspace
$$
\langle  d \bm p d \bm x \otimes  x_1^3,  d \bm p d \bm x \otimes x_2^3 \rangle.
$$
Following the ideas already discussed, the output of the mirror map should carry the variable $p_2$, in order to record the action on the $X_i$ variables. Furthermore, it is natural to keep track of the two distinct generators $d \bm p d \bm x \otimes  x_1^3$ and $d \bm p d \bm x \otimes x_2^3$ via their polynomial part. Therefore, we define
\begin{equation} \label{dim 2 case 1}
d \bm p d \bm x \otimes  x_i^3\left|\frac{1}{9}(0,0,0,3,3,3;0,0)\right\rangle \mapsto d \bm p d \bm x d \bm X\otimes p_2 x_i^3|id\rangle .
\end{equation}
Notice that this assignment is valid also for $i=3$.

Similarly, we define
\begin{equation} \label{dim 2 case 2}
d \bm p d \bm X \otimes  X_i^3\left|\frac{1}{9}(0,0,0,3,3,3;0,0)\right\rangle \mapsto d \bm p d \bm x d \bm X\otimes p_2 X_i^3|id\rangle
\end{equation}
Notice that the assignments in equations \eqref{dim 2 case 1} and \eqref{dim 2 case 2} could be regarded as following an extension of the rule in equation \eqref{assignment rule}, where we allow some degree of complexity in the polynomial part of the input as well.

Thus, we are left with assigning a mirror element for the three elements arising in from $\gamma \in \Gamma^T$, with $r_\gamma=2$ and $n_\gamma =0$. In the case of the elements
$$dt\left|\frac{1}{9}(3,3,3,6,6,6;0,0)\right\rangle, \quad dt\left|\frac{1}{9}(6,6,6,3,3,3;0,0)\right\rangle $$
we can follow the algorithm in equation \eqref{assignment rule}. In this way, we get
$$
dt\left|\frac{1}{9}(3,3,3,6,6,6;0,0)\right\rangle \mapsto d \bm p d \bm x d \bm X \otimes p_2 X_1 X_2 X_3 | id \rangle,
$$
and
$$
dt\left|\frac{1}{9}(6,6,6,3,3,3;0,0)\right\rangle \mapsto d \bm p d \bm x d \bm X \otimes p_1 x_1 x_2 x_3| id \rangle.
$$

Now, we have to consider the element
$$t dt\left|\frac{1}{9}(3,3,3,3,3,3;0,0)\right\rangle.$$
This element does not follow directly in the patter of equation \eqref{assignment rule}. On the other hand, it comes from the most symmetric element $\gamma \in \Gamma^T \setminus \lbrace 0 \rbrace$. Analogously, there is a highly symmetric generator of $\mathcal{H}^{2,1}_\Gamma(W_1,W_2)$, namely
$$
d \bm p d \bm x d \bm X\otimes p_1 X_1 X_2 X_3 = d \bm p d \bm x d \bm X\otimes p_2 x_1 x_2 x_3 = d \bm p d \bm x d \bm X \otimes p_2 X_i^3 = d \bm p d \bm x d \bm X\otimes p_1 x_i^3.
$$
Given the common symmetries, we suggest the choice
$$
t dt\left|\frac{1}{9}(3,3,3,3,3,3;0,0)\right\rangle \mapsto d \bm p d \bm x d \bm X\otimes p_1 X_1 X_2 X_3 | id \rangle. 
$$

As the above assignments determine a bijective correspondence between a basis of $\mathcal{H}^{1,1}_{\Gamma^T}(W_1,W_2)$ and one of $\mathcal{H}^{2,1}_{\Gamma}(W_1,W_2)$, we can extend the map by linearity. This completes the definition of the mirror map in Definition \ref{main result}.

\appendix

\section{Tables}

Here we collect the tables referenced above.
\begin{center}
\rowcolors{2}{light-gray}{white}
\begin{table}[H]
\begin{tabular}{ l  c } 
\toprule
Generator & \# \\
\midrule
  $p_1X_1X_2X_3=p_2x_1x_2x_3=p_2X_i^3=p_1x_i^3$ & 1  \\ 
  $p_1x_1x_2x_3 = p_1 \sum X_i^3$ & 1 \\ 
	$p_2X_1X_2X_3 = p_2 \sum x_i^3$ & 1\\
$p_1X_1^3$  & 1 \\
$p_1X_2^3$ & 1\\
$p_2x_1^3$ & 1\\
$p_2x_2^3$ & 1\\
$p_1x_ix_jX_k$ & 9\\
$p_2x_iX_jX_k$ & 9\\
$p_1X_iX_jx_k  =p_2X_l^2x_k$ & 9\\
$p_1x_i^2X_j  =p_2x_jx_kX_j$ & 9\\
$p_1X_i^2x_j$ & 9 \\
$p_2x_i^2X_j$ & 9 \\
$p_1x_i^2x_j =p_2 x_j^2x_k$ & 6\\
$p_1X_i^2X_j=p_2X_iX_k^2$ & 6 \\
\bottomrule\addlinespace
\end{tabular}
\captionsetup{width =0.7\LTcapwidth}
\caption{A basis of $H_{\Gamma}^{2,1}(W_1,W_2)$. All indices $i,j,k,l$ range from 1 to 3, distinct indices have distinct values. On every row there are all the possible ways of writing a generator, and the number of generators of that form.}
\label{tab_2,1}
\end{table}
\end{center}

\rowcolors{2}{light-gray}{white}
\begin{table}[H]
\begin{minipage}[b]{0.45\linewidth}\centering
\begin{tabular}{ cccc } 
\toprule
$r_\gamma$ & $n_\gamma$ & element & $a_\gamma$ \\
\midrule
  2 & 6 & (0,...,0) & 0\\
  \midrule
  2 & 3 & $\frac19(0,0,0,3,3,3;0,0)$ & 1\\
  2 & 3 & $\frac19(3,3,3,0,0,0;0,0)$ & 1\\
  2 & 3 & $\frac19(0,0,0,6,6,6;0,0)$ & 2\\
  2 & 3 & $\frac19(0,0,0,6,6,6;0,0)$ & 2\\
  
\bottomrule\addlinespace
\end{tabular}
\captionsetup{width =0.6\LTcapwidth}
\caption{Elements of $\Gamma^T$ for which $r_\gamma < n_\gamma$}
\label{tab_JacobiGamma^vee}
\end{minipage}
\hspace{0.5cm}
\rowcolors{2}{light-gray}{white}
\begin{minipage}[b]{0.45\linewidth}
\centering
\begin{tabular}{ ccccc } 
\toprule
$r_\gamma$ & $n_\gamma$ & element type & $a_\gamma$ & \# \\

  \midrule
  2 & 0 & $\frac19(3,3,3,3,3,3;0,0)$ & 2 & 1\\
  2 & 0 & $\frac19(6,6,6,3,3,3;0,0)$ & 3 & 1\\
  2 & 0 & $\frac19(3,3,3,6,6,6;0,0)$ & 3 & 1\\
  2 & 0 & $\frac19(6,6,6,6,6,6;0,0)$ & 4 & 1\\
  \midrule
  2 & 1 & $\frac19(3,3,3,3,0,6;0,0)$ & 2 & 6\\
  2 & 1 & $\frac19(3,0,6,3,3,3;0,0)$ & 2 & 6\\
  2 & 1 & $\frac19(6,6,6,3,0,6;0,0)$ & 3 & 6\\
  2 & 1 & $\frac19(3,0,6,6,6,6;0,0)$ & 3 & 6\\
  \midrule
  1 & 0 & $\frac19(2,2,5,6,3,6;3,0)$ & 3 & 9\\
  1 & 0 & $\frac19(6,3,6,2,2,5;0,3)$ & 3 & 9\\
  1 & 0 & $\frac19(3,6,3,1,7,1;0,6)$ & 3 & 9\\
  1 & 0 & $\frac19(1,7,1,3,6,3;6,0)$ & 3 & 9\\
  1 & 0 & $\frac19(3,6,3,4,4,1;0,6)$ & 3 & 9\\
  1 & 0 & $\frac19(4,4,1,3,6,3;6,0)$ & 3 & 9\\
  
  1 & 0 & $\frac19(6,6,3,8,2,8;0,3)$ & 4 & 9\\
  1 & 0 & $\frac19(8,2,8,6,6,3;3,0)$ & 4 & 9\\
  1 & 0 & $\frac19(6,6,3,5,5,8;0,3)$ & 4 & 9\\
  1 & 0 & $\frac19(5,5,8,6,6,3;3,0)$ & 4 & 9\\
  1 & 0 & $\frac19(7,4,7,3,3,6;6,0)$ & 4 & 9\\
  1 & 0 & $\frac19(3,3,6,7,4,7;0,6)$ & 4 & 9\\
  
\bottomrule\addlinespace
\end{tabular}
\captionsetup{width =0.6\LTcapwidth}
\caption{Elements of $\Gamma^T$ for which $r_\gamma \geq n_\gamma$.}
\label{tab_projGamma^vee}
\end{minipage}
\end{table}

\rowcolors{2}{light-gray}{white}
\begin{longtable}[H]{ cccc}
\rowcolor{white}
 \multicolumn{1}{l}{\textsc{Table \ref*{tab_complete_mirrormap}}} & \\
\toprule \rowcolor{white}
$H_{\Gamma^T}^{1,1}(W_1,W_2)$ & $\rightarrow$ &  $H_{\Gamma}^{2,1}(W_1,W_2)$ \\
\midrule\rowcolor{white}
$tdt|\frac{1}{9}(3,3,3,3,3,3;0,0)\rangle$ & $\mapsto$ & $p_1x_1^3|id\rangle$\\
\midrule
$dt|\frac{1}{9}(3,3,3,6,6,6;0,0)\rangle$ & $\mapsto$ & $p_2X_1X_2X_3|id\rangle$\\
$dt|\frac{1}{9}(6,6,6,3,3,3;0,0)\rangle$ & $\mapsto$ & $p_1x_1x_2x_3|id\rangle$\\
\midrule
$d\bm p d\bm x x_1^3|\frac{1}{9}(0,0,0,3,3,3;0,0)\rangle$&$ \mapsto$&$  p_2x_1^3|id\rangle$ \\
$d\bm p d\bm x x_2^3|\frac{1}{9}(0,0,0,3,3,3;0,0)\rangle$&$ \mapsto $&$ p_2x_2^3|id\rangle$\\
$d\bm pd\bm X X_1^3|\frac{1}{9}(3,3,3,0,0,0;0,0)\rangle$&$ \mapsto$&$  p_1X_1^3|id\rangle$\\
$d\bm p d\bm X X_2^3|\frac{1}{9}(3,3,3,0,0,0;0,0)\rangle$&$ \mapsto$&$  p_1X_2^3|id\rangle$\\
\midrule
$dt|\frac19(3,3,3,3,0,6;0,0)\rangle$&$ \mapsto$&$ p_1X_1X_3^2|id\rangle$ \\
$dt| \frac19(3, 3, 3, 0, 3, 6; 0,0) \rangle $&$ \mapsto  $&$p_1X_2X_3^2|id\rangle$\\
$dt| \frac19(3, 3, 3, 0, 6, 3; 0,0) \rangle $&$ \mapsto  $&$p_1X_2^2X_3|id\rangle$\\
$dt| \frac19(3, 3, 3, 6, 0, 3; 0,0) \rangle $&$ \mapsto  $&$p_1X_1^2X_3|id\rangle$\\
$dt| \frac19(3, 3, 3, 6, 3, 0; 0,0) \rangle $&$ \mapsto  $&$p_1X_1^2X_2|id\rangle$\\
$dt| \frac19(3, 3, 3, 3, 6, 0; 0,0) \rangle $&$ \mapsto  $&$p_1X_1X_2^2|id\rangle$\\
\midrule
$dt| \frac19(0, 3, 6, 3, 3, 3; 0,0) \rangle $&$ \mapsto  $&$p_2X_2X_3^2|id\rangle$\\
$dt| \frac19(3, 0, 6, 3, 3, 3; 0,0) \rangle $&$ \mapsto  $&$p_2X_1X_3^2|id\rangle$\\
$dt| \frac19(6, 0, 3, 3, 3, 3; 0,0) \rangle $&$ \mapsto  $&$p_2X_1^2X_3|id\rangle$\\
$dt| \frac19(0, 6, 3, 3, 3, 3; 0,0) \rangle $&$ \mapsto  $&$p_2X_2^2X_3|id\rangle$\\
$dt| \frac19(3, 6, 0, 3, 3, 3; 0,0) \rangle $&$ \mapsto  $&$p_2X_1X_2^2|id\rangle$\\
$dt| \frac19(6, 3, 0, 3, 3, 3; 0,0) \rangle $&$ \mapsto  $&$p_2X_1^2X_2|id\rangle$\\
\midrule
$dt|\frac19(2,2,5,6,3,6;3,0)\rangle$&$\mapsto$&$ p_2x_3X_1X_3|id\rangle$\\
$dt| \frac19(2, 2, 5, 3, 6, 6; 3,0  ) \rangle $&$ \mapsto $&$ p_2x_3X_2X_3|id\rangle $\\
$dt| \frac19(2, 2, 5, 6, 6, 3; 3,0 ) \rangle $&$ \mapsto $&$ p_2x_3X_1X_2|id\rangle $\\
$dt| \frac19(2, 5, 2, 6, 3, 6; 3,0 ) \rangle $&$ \mapsto $&$ p_2x_2X_1X_3|id\rangle $\\
$dt| \frac19(2, 5, 2, 3, 6, 6; 3,0 ) \rangle $&$ \mapsto $&$ p_2x_2X_2X_3|id\rangle $\\
$dt| \frac19(5, 2, 2, 6, 3, 6; 3,0 ) \rangle $&$ \mapsto $&$ p_2x_1X_1X_3|id\rangle$\\
$dt| \frac19(5, 2, 2, 3, 6, 6; 3,0 ) \rangle $&$ \mapsto $&$ p_2x_1X_2X_3|id\rangle$\\
$dt| \frac19(2, 5, 2, 6, 6, 3; 3,0 ) \rangle $&$ \mapsto $&$ p_2x_2X_1X_2|id\rangle $\\
$dt| \frac19(5, 2, 2, 6, 6, 3; 3,0 ) \rangle $&$ \mapsto $&$ p_2x_1X_1X_2|id\rangle $\\
\midrule
$dt| \frac19( 6, 3, 6, 2, 2, 5; 0,3 ) \rangle $&$ \mapsto $&$ p_1x_1x_3X_3|id\rangle$\\
$dt| \frac19(3, 6, 6, 2, 2, 5; 0,3 ) \rangle $&$ \mapsto $&$p_1x_2x_3X_3|id\rangle $\\
$dt| \frac19( 6, 3, 6, 2, 5, 2; 0,3 ) \rangle $&$ \mapsto $&$p_1x_1x_3X_2|id\rangle $\\
$dt| \frac19(6, 3, 6, 5, 2, 2; 0,3 )\rangle $&$ \mapsto $&$ p_1x_1x_3X_1|id\rangle$\\
$dt| \frac19(3, 6, 6, 2, 5, 2; 0,3 ) \rangle $&$ \mapsto $&$p_1x_2x_3X_2|id\rangle $\\
$dt| \frac19(3, 6, 6, 5, 2, 2; 0,3 ) \rangle $&$ \mapsto $&$p_1x_2x_3X_1|id\rangle $\\
$dt| \frac19(6, 6, 3, 2, 2, 5; 0,3 ) \rangle $&$ \mapsto $&$p_1x_1x_2X_3|id\rangle $\\
$dt| \frac19(6, 6, 3, 2, 5, 2; 0,3 ) \rangle $&$ \mapsto $&$p_1x_1x_2X_2|id\rangle $\\
$dt| \frac19(6, 6, 3, 5, 2, 2; 0,3 ) \rangle $&$ \mapsto $&$p_1x_1x_2X_1|id\rangle $\\
\midrule
$dt| \frac19(3, 6, 3, 1, 7, 1; 0,6 ) \rangle $&$ \mapsto $&$ p_1x_2X_2^2id\rangle$\\
$dt| \frac19(3, 6, 3, 7, 1, 1; 0,6 ) \rangle $&$ \mapsto $&$ p_1x_2X_1^2id\rangle$\\
$dt| \frac19(6, 3, 3, 1, 7, 1; 0,6 ) \rangle $&$ \mapsto $&$ p_1x_1X_2^2id\rangle$\\
$dt| \frac19(6, 3, 3, 7, 1, 1; 0,6 ) \rangle $&$ \mapsto $&$ p_1x_1X_1^2id\rangle$\\
$dt| \frac19(3, 3, 6, 1, 1, 7; 0,6 ) \rangle $&$ \mapsto $&$ p_1x_3X_3^2id\rangle$\\
$dt| \frac19(3, 6, 3, 1, 1, 7; 0,6 ) \rangle $&$ \mapsto $&$ p_1x_2X_3^2id\rangle$\\
$dt| \frac19(6, 3, 3, 1, 1, 7; 0,6 ) \rangle $&$ \mapsto $&$ p_1x_1X_3^2id\rangle$\\
$dt| \frac19(3, 3, 6, 1, 7, 1; 0,6 ) \rangle $&$ \mapsto $&$ p_1x_3X_2^2id\rangle$\\
$dt| \frac19(3, 3, 6, 7, 1, 1; 0,6 ) \rangle $&$ \mapsto $&$ p_1x_3X_1^2id\rangle$\\
\midrule
$dt| \frac19(3, 6, 3, 4, 4, 1; 0,6 ) \rangle $&$ \mapsto $&$ p_1x_2X_1X_2id\rangle $\\
$dt| \frac19(6, 3, 3, 4, 4, 1; 0,6 ) \rangle $&$ \mapsto $&$p_1x_1X_1X_2id\rangle $\\
$dt| \frac19(3, 3, 6, 4, 1, 4; 0,6 ) \rangle $&$ \mapsto $&$p_1x_3X_1X_3id\rangle $\\
$dt| \frac19(3, 3, 6, 1, 4, 4; 0,6 ) \rangle $&$ \mapsto $&$p_1x_3X_2X_3id\rangle $\\
$dt| \frac19(3, 6, 3, 4, 1, 4; 0,6 ) \rangle $&$ \mapsto $&$p_1x_2X_1X_2id\rangle $\\
$dt| \frac19(3, 6, 3, 1, 4, 4; 0,6 ) \rangle $&$ \mapsto $&$p_1x_2X_2X_3id\rangle $\\
$dt| \frac19(6, 3, 3, 4, 1, 4; 0,6 ) \rangle $&$ \mapsto $&$p_1x_1X_1X_3id\rangle $\\
$dt| \frac19(6, 3, 3, 1, 4, 4; 0,6 ) \rangle $&$ \mapsto $&$p_1x_1X_2X_3id\rangle $\\
$dt| \frac19(3, 3, 6, 4, 4, 1; 0,6 ) \rangle $&$ \mapsto $&$p_1x_3X_1X_2id\rangle $\\
\midrule
$dt| \frac19( 1, 1, 7, 3, 3, 6; 6,0 ) \rangle $&$ \mapsto $&$p_2x_3^2X_3id\rangle $\\
$dt| \frac19( 1, 7, 1, 3, 6, 3; 6,0 ) \rangle $&$ \mapsto $&$p_2x_2^2X_2id\rangle $\\
$dt| \frac19(1, 7, 1, 6, 3, 3; 6,0 ) \rangle $&$ \mapsto $&$p_2x_2^2X_1id\rangle $\\
$dt| \frac19(7, 1, 1, 3, 6, 3; 6,0 ) \rangle $&$ \mapsto $&$p_2x_1^2X_2id\rangle $\\
$dt| \frac19(7, 1, 1, 6, 3, 3; 6,0 ) \rangle $&$ \mapsto $&$p_2x_1^2X_1id\rangle $\\
$dt| \frac19(1, 1, 7, 3, 6, 3; 6,0 ) \rangle $&$ \mapsto $&$p_2x_3^2X_2id\rangle $\\
$dt| \frac19(1, 1, 7, 6, 3, 3; 6,0 ) \rangle $&$ \mapsto $&$p_2x_3^2X_1id\rangle $\\
$dt| \frac19(1, 7, 1, 3, 3, 6; 6,0 ) \rangle $&$ \mapsto $&$p_2x_2^2X_3id\rangle $\\
$dt| \frac19(7, 1, 1, 3, 3, 6; 6,0 ) \rangle $&$ \mapsto $&$p_2x_1^2X_3id\rangle$\\
\midrule
$dt| \frac19(4, 4, 1, 3, 6, 3; 6,0 ) \rangle $&$ \mapsto $&$p_2x_1x_2X_2id\rangle $\\
$dt| \frac19( 4, 4, 1, 6, 3, 3;6,0   ) \rangle $&$ \mapsto $&$p_2x_1x_2X_1id\rangle $\\
$dt| \frac19(4, 1, 4, 3, 3, 6;6,0   ) \rangle $&$ \mapsto $&$ p_2x_1x_3X_3id\rangle$\\
$dt| \frac19(1, 4, 4, 3, 3, 6;6,0   ) \rangle $&$ \mapsto $&$p_2x_2x_1X_3id\rangle$\\
$dt| \frac19(4, 4, 1, 3, 3, 6;6,0  ) \rangle $&$ \mapsto $&$ p_2x_1x_2X_3id\rangle$\\
$dt| \frac19( 4, 1, 4, 3, 6, 3; 6,0  ) \rangle $&$ \mapsto $&$p_2x_1x_3X_2id\rangle $\\
$dt| \frac19( 4, 1, 4, 6, 3, 3; 6,0 ) \rangle $&$ \mapsto $&$ p_2x_1x_3X_1id\rangle$\\
$dt| \frac19(1, 4, 4, 3, 6, 3;6,0   ) \rangle $&$ \mapsto $&$p_2x_2x_3X_2id\rangle $\\
$dt| \frac19(1, 4, 4, 6, 3, 3;6,0   ) \rangle $&$ \mapsto $&$p_2x_2x_3X_1id\rangle $\\
\bottomrule \addlinespace
\captionsetup{width =0.7\LTcapwidth}\rowcolor{white}
\caption{The complete description of the mirror map on the generators of $\mathcal{H}^{1,1}_{\Gamma^T}$, grouped according to the subdivision described above}
\label{tab_complete_mirrormap}
\end{longtable}

\newpage

\printbibliography

\Addresses

\end{document}